\documentclass[manuscript,screen,nonacm]{acmart}
\usepackage{amsmath,bm}
\usepackage{listings}
\usepackage{subfig}
\usepackage{multicol}

\begin{document}

\title[FdeSolver: A Julia Package for Solving FDEs]{FdeSolver: A Julia Package for Solving Fractional Differential Equations}

\author{Moein Khalighi}
\email{moein.khalighi@utu.fi}
\orcid{0000-0001-8176-0367}
\affiliation{%
  \institution{Department of Computing, Faculty of Technology, University of Turku}
  \streetaddress{Yliopistonmäki}
  \city{Turku}
  \country{Finland}
  \postcode{FI-20014}
}

\author{Giulio Benedetti}
\orcid{0000-0002-8732-7692}
\affiliation{%
\institution{Department of Computing, Faculty of Technology, University of Turku}
  \streetaddress{Yliopistonmäki}
  \city{Turku}
  \country{Finland}}

\author{Leo Lahti}
\orcid{0000-0001-5537-637X}
\affiliation{%
  \institution{Department of Computing, Faculty of Technology, University of Turku}
  \city{Turku}
  \country{Finland}
}

\renewcommand{\shortauthors}{Khalighi, et al.}

\begin{abstract}
Implementing and executing numerical algorithms to solve fractional differential equations has been less straightforward than using their integer-order counterparts, posing challenges for practitioners who wish to incorporate fractional calculus in applied case studies. Hence, we created an open-source Julia package, FdeSolver, that provides numerical solutions for fractional-order differential equations based on product-integration rules, predictor–corrector algorithms, and the Newton-Raphson method. The package covers solutions for one-dimensional equations with orders of positive real numbers. For high-dimensional systems, the orders of positive real numbers are limited to less than (and equal to) one. Incommensurate derivatives are allowed and defined in the Caputo sense. Here, we summarize the implementation for a representative class of problems, provide comparisons with available alternatives in Julia and Matlab, describe our adherence to good practices in open research software development, and demonstrate the practical performance of the methods in two applications; we show how to simulate microbial community dynamics and model the spread of Covid-19 by fitting the order of derivatives based on epidemiological observations. Overall, these results highlight the efficiency, reliability, and practicality of the FdeSolver Julia package.
\end{abstract}



\maketitle

\section{Introduction}

Fractional calculus is the study of non-integer order derivatives and integrals. Despite its origins in pure mathematics, fractional derivatives have recently been applied to a wide range of real-world scenarios~\cite{de2014review}. Furthermore, advanced computational approaches and new generations of computers have brought fractional calculus to a diverse spectrum of disciplines~\cite{SUN-reviewApplication}. Because the domain of application of fractional calculus has expanded, but analytical solutions are sometimes difficult or even impossible to find, researchers are studying, developing, and using numerical methods to solve fractional differential equations (FDEs).

One significant application of fractional derivatives is the incorporation of memory effects into the classical system with integer order~\cite{Saeedian2017PhysRevE, Safdari_plosone, Eftekhari2022, KhalighiPlosCB,KhalighiSymmetry,AMIRIAN2020}. Although memory effects in real phenomena do not have a unique shape, their common feature is the incorporation of an average of the previous values into the current value. By holding this feature, the concept of memory effects has been interpreted in different ways by the various definitions of fractional derivatives~\cite{loverro2004fractional}. Caputo has defined one of the most popular differential operators~\cite{podlubny1998fractional,kilbas2006theory}, as it has two specific advantages: 1) the derivative of a constant is zero, and 2) it includes convolution integrals with a singular power-law kernel; these are desirable properties in applied mathematics and other practical areas. The former leads to having an operator for modelling systems with higher orders and nonzero initial values, and the latter is applicable to many real-world phenomena with gradually decaying memory effects~\cite{Diethelm2020Nonsingular}.

FdeSolver is one of the first open-source Julia packages for solving FDEs. It can solve models of ordinary differential equations (ODEs) with Caputo fractional derivatives. FractionalDiffEq~\cite{SciFracx} is an alternative Julia package that is currently under development and provides solutions to FDEs with some similarities to our methodology. Our benchmarking experiments compare the methods shared by these two implementations and demonstrate that a comparatively better performance is achieved with the FdeSolver package.

We have developed the FdeSolver package based on well-established mathematical algorithms~\cite{diethelm2002predictor,diethelm2004detailed,garrappa2010linear} that were formerly designed and implemented as Matlab routines~\cite{Garrappa2018, Matlab-Review-FDEs}. These implementations convert the FDE problem to a Volterra integral equation, discretize it via product-integration (PI) rules, numerically solve this by using either the predictor-corrector (PC) method or the Newton–Raphson (NR) method, and finally, reduce the cost of computations by fast Fourier transform (FFT).  However, because Matlab is a proprietary programming language, its usability is restricted to the holders of a Matlab license. Here, we provide an open and independent alternative to the Matlab routines, which depend on libraries written in C, C++, or Fortran. 

In addition to showing similar or improved levels of performance in our numerical experiments when compared to the analogous implementations in Matlab and Julia, the FdeSolver package has additional benefits, as Julia is a generic open-source programming language, it has a built-in package manager, and its multiple dispatch paradigm paves the way for interoperability with other packages within the larger Julia ecosystem. We believe that these features can make our implementation attractive for a broad range of researchers~\cite{Bezanson_Julia}, and support further package development, for instance, to reduce the computation time to the standards of C and C++.

We organized the paper as follows: We start with preliminary mathematics and numerical methods in Sec. \ref{Sec:PC and NR}. Then, we present the structure of FdeSolver in Sec. \ref{Sec: software}. The usage of the package for five classes of problems is demonstrated in Sec. \ref{Sec: usage}, including 
benchmarks with alternative solvers. We provide two applications of the package in the context of ecological and epidemiological models in Sec. \ref{sec: applications}. We conclude with Sec. \ref{sec: Conclusion}.

\section{Preliminaries and Numerical Scheme}\label{Sec:PC and NR}
This section focuses on some primary definitions in fractional calculus based on Refs.~\cite{podlubny1998fractional, kilbas2006theory} and the numerical methods used in the package, which is proposed by Diethelm~\cite{diethelm2002predictor, diethelm2004detailed} and implemented in Matlab by Garrappa~\cite{Garrappa2018}. 

Let us suppose \(\mathcal{D}_{t_0}^{\beta}\) as a notation of time-fractional Caputo derivative with order of \(\beta \in \mathbb{R}^+\) from initial value \(t_0\). Thus, the fractional derivative of a given differentiable function \(g\) at time \(t\) in the sense of Caputo is defined as:

\begin{equation}\label{eq:caputo}
\mathcal{D}_{t_0}^{\beta} g(t) = I^{1-\beta}_{t_0}g'(t)=\frac{1}{\Gamma (1- \beta  )}\int_{{t_0}}^{t}{\frac{g'(\tau )d\tau }{{(t - {\tau})}^{\beta}}},
\end{equation}
in which \(I^{1-\beta}_{t_0}\) is Riemann-Liouville fractional integral of order \(1-\beta\) that is defined by
\begin{equation}\label{eq:RLInt}
I^{\beta}_{t_0}g(t)=\frac{1}{\Gamma (\beta )}\int_{{t_0}}^{t}{\frac{g(\tau )d\tau }{{(t - {\tau})}^{1-\beta}}},
\end{equation}
where $\Gamma$ denotes the gamma function. 

FdeSolver can solve the following types of equations and systems

\begin{equation}\label{eq:scalar}
    \mathcal{D}_{t_0}^{\beta_i} X_i = f_i(t, \bm{X}), \quad i=1,...,M,
\end{equation}
with the initial conditions \(X_i(t_0)=X_{i,0}\) when \(0<\beta_i \leq 1\) for a system of equations, and \(X_1^{(j)}(t_0)=X^{(j)}_{1,0}\) \((j=0,...,m-1)\) for one-dimensional equations when \(\beta_1>1\) and \(m\) is the smallest integer greater than or equal to the order derivative,~\(\beta_1\). We can rewrite the fractional order system \eqref{eq:scalar} in the following vector form

\begin{equation}\label{eq:matrix}
    \mathcal{D}_{t_0}^{\bm{\beta}}\bm{X}=\bm{F}(t,\bm{X}),
\end{equation}
where \(\bm{X}=[X_1,...,X_M]\), \(\bm{\beta}=[\beta_1,...,\beta_M]\), and \(\bm{F}=[f_1,...,f_M]\).

The initial value problem \eqref{eq:matrix} is equivalent to the Volterra integral equation~\cite{kilbas2006theory,diethelm2002predictor}
\begin{equation}\label{eq:VIE}
    \bm{X}(t)=T_{m-1}[\bm{X};t_0]+\frac{1}{\Gamma(\bm{\beta})}\int_{t_0}^t(t-\tau)^{\bm{\beta}-1}\bm{F}(\tau,\bm{X}(\tau))d\tau.
\end{equation}

in which \(T_{m-1}[X_i;t_0]\) is Taylor polynomial of degree \(m-1\) for the function \(X_i(t)\) centered at \(t_0\) defined as

\begin{displaymath}
 T_{m-1}[X_i;t_0](t_n)=\sum_{k=0}^{m-1}\frac{(t-t_0)^{k}}{k!}X^{(k)}_i(t_0).
\end{displaymath}

The integrals in Equations \eqref{eq:VIE} give us history-dependent dynamics, and the presence of power in the kernels provides a nonlocal feature of our fractional order modellings. These, however, cause a non-smooth behaviour at the initial time, thus straightforward numerical methods, such as polynomial approximations, cannot achieve the solution, and only those are acceptable that each step of computation involves the whole history of the solution.

Therefore, based on product-integration rules~\cite{diethelm2002predictor, garrappa2010linear}, we discretise the integral term of Eq.~\eqref{eq:VIE} in \(n\) points and the step size \(h>0\); \(t_r=t_0+rh\) \((r=0,...,n-1)\), so that the Eq.~\eqref{eq:VIE} turns to

\begin{equation}\label{eq:disVIE}
    \bm{X}(t_n)=T_{m-1}[\bm{X};t_0](t_n)+\frac{1}{\Gamma(\bm{\beta})}\sum_{r=0}^{n-1}\int_{t_r}^{t_{r+1}}(t_n-\tau)^{\bm{\beta}-1}\bm{F}(\tau,\bm{X}(\tau))d\tau.
\end{equation}
Then we use piece-wise interpolating polynomials for the approximation of integral terms.

\subsection{Predictor-Corrector Method}\label{Sec: pc}
The predictor-corrector approach has been proposed by using a generalization of Adams multi-step methods for the solution of fractional ordinary differential equations~\cite{diethelm2002predictor, garrappa2010linear}. To use this technique for the solution of Eq.~\eqref{eq:disVIE}, we start with a preliminary explicit estimation as a {\itshape predictor} and then improve it by the iteration of implicit approximations as a {\itshape corrector}. Hence, the PI rectangular rule~\cite{Garrappa2018} gives the matrix form of the predictor of the solution of Equation \eqref{eq:disVIE}:
\begin{equation}\label{explicit}
    \bm{X}^{[0]}=\bm{X}_0+\bm{h}\bm{b}\bm{F}^{n-1}_{0},
\end{equation}
and the PI trapezoidal rule~\cite{Garrappa2018} provides the corrector:
\begin{equation}\label{implicit}
    \bm{X}^{[\mu]}=\bm{X}_0+\bm{h}\bm{c}\bm{F}_{0}+\bm{h}\bm{d}\bm{F}_{1}^{n},
\end{equation}
in details, Eq \eqref{explicit} is
\begin{displaymath}
 \begin{bmatrix}
    X_1^{[0]}(t_n) \\
    X_2^{[0]}(t_n) \\
    \vdots \\
    X_M^{[0]}(t_n)
\end{bmatrix}^{\intercal}
=
\begin{bmatrix}
    T(X_1) \\
    T(X_2)\\
    \vdots \\
    T(X_M)
\end{bmatrix}  ^{\intercal}
+
\begin{bmatrix}
    h^{\beta_1} \\
    h^{\beta_2} \\
    \vdots \\
    h^{\beta_M}
\end{bmatrix}^{\intercal}
\begin{bmatrix}
  b_{1,n-1} & b_{1,n-2} & \dots  & b_{1,0} \\
  b_{2,n-1} & b_{2,n-2} & \dots  & b_{2,0} \\
  \vdots &\vdots & \ddots & \vdots \\
  b_{M,n-1} & b_{M,n-2} & \dots  & b_{M,0} \\
\end{bmatrix}  
\begin{bmatrix}
    f_1(\bm{X}_0) & \dots  & f_M(\bm{X}_0) \\        f_1(\bm{X}_1) & \dots  & f_M(\bm{X}_1) \\
    \vdots & \ddots & \vdots \\
    f_1(\bm{X}_{n-1}) & \dots  & f_M(\bm{X}_{n-1}) \\
\end{bmatrix}  
\end{displaymath}
where 
\begin{math}
 T(X_i)=T_{m-1}[X_i;t_0](t_n)
 \end{math}
 and
 \begin{math}
 b_{i,r}=\frac{(r+1)^{\beta_i}-(r)^{\bm{\beta}_i}}{\Gamma(\bm{\beta}_i+1)},\; i=1,...,M,\; r=0,...,n-1,
 \end{math}
and Eq.~\eqref{implicit} is
\begin{align*}
 \begin{bmatrix}
    X_1^{[\mu]}(t_n) \\
    X_2^{[\mu]}(t_n) \\
    \vdots \\
    X_M^{[\mu]}(t_n)
\end{bmatrix}^{\intercal}
&=                    
\begin{bmatrix}
    T(X_1) \\
    T(X_2)\\
    \vdots \\
    T(X_M)
\end{bmatrix}  ^{\intercal}
+
\begin{bmatrix}
    h^{\beta_1}\times c_{1,n} \times f_1(\bm{X}_0) \\
    h^{\beta_2}\times c_{2,n} \times f_2(\bm{X}_0) \\
    \vdots \\
    h^{\beta_M} \times c_{M,n} \times f_M(\bm{X}_0)
\end{bmatrix}^{\intercal} \notag \\ 
&\quad +
\begin{bmatrix}
    h^{\beta_1} \\
    h^{\beta_2} \\
    \vdots \\
    h^{\beta_M}
\end{bmatrix}^{\intercal}
\begin{bmatrix}
  d_{1,n-1} & d_{1,n-2} & \dots  & d_{1,0} \\
  d_{2,n-1} & d_{2,n-2} & \dots  & d_{2,0} \\
  \vdots &\vdots & \ddots & \vdots \\
  d_{M,n-1} & d_{M,n-2} & \dots  & d_{M,0} \\
\end{bmatrix}  
 \begin{bmatrix}
    f_1(\bm{X}_1) & \dots  & f_M(\bm{X}_1) \\
    \vdots & \ddots & \vdots \\
    f_1(\bm{X}_{n-1}) & \dots  & f_M(\bm{X}_{n-1}) \\        f_1(\bm{X}^{[\mu-1]}_n) & \dots  & f_M(\bm{X}^{[\mu-1]}_n) \\
\end{bmatrix}
,\; \mu=1,2,...
\end{align*}
in which
\begin{displaymath}
c_{i,n}=\frac{(n-1)^{\beta_i+1}-n^{\beta_i}(n-\beta_i-1)}{\Gamma(\beta_i+2)},
\end{displaymath}
\begin{displaymath}
 d_{i,r}= 
\begin{cases}
 \frac{1}{\Gamma({{\beta_i}+2})},& \text{if } r=0,\\
 \frac{(r-1)^{{\beta_i}+1}-2r^{{\beta_i}+1}+(r+1)^{{\beta_i}+1}}{\Gamma({\beta_i}+2)},              & \text{if}  \;r=1,2, ..., n-1.
\end{cases}
\end{displaymath}
and, \(\mu\) is the number of corrections. This is the main algorithm of our package with the convergence rate of \(\mathcal{O}(h)\) for the predictor and \(\mathcal{O}(h^{min\{1+\beta,2\}})\) for the corrector~\cite{Garrappa2018}. 

\subsection{Newton–Raphson Method}\label{Sec: NR}
The predictor-corrector method may not be sufficient for stiff problems, which need too small step size \(h\) leading to large computation costs for the sake of accuracy. Hence, we use the more efficient approximation with better stability properties, that is the implicit method Eq.~\eqref{implicit}. But, this method needs a prior approximation of the current step besides the computed values of the previous steps. So, we rewrite Eq.~\eqref{implicit} as:

\begin{equation}\label{eq: IM}
    \bm{X}_n=\bm{\Psi}_{n-1}+\bm{a}_0\bm{F}_n
\end{equation}
where \(\bm{\Psi}_{n-1}\) denotes the term including all the explicitly known information of the previous steps
\begin{displaymath}
\bm{\Psi}_{n-1}
=  
\begin{bmatrix}
    T(X_1) \\
    T(X_2)\\
    \vdots \\
    T(X_M)
\end{bmatrix}  ^{\intercal}
+
\begin{bmatrix}
    h^{\beta_1}\times c_{1,n} \times f_1(\bm{X}_0) \\
    h^{\beta_2}\times c_{2,n} \times f_2(\bm{X}_0) \\
    \vdots \\
    h^{\beta_M} \times c_{M,n} \times f_M(\bm{X}_0)
\end{bmatrix}^{\intercal} 
+
\begin{bmatrix}
    h^{\beta_1} \\
    h^{\beta_2} \\
    \vdots \\
    h^{\beta_M}
\end{bmatrix}^{\intercal}
\begin{bmatrix}
  d_{1,n-1} & d_{1,n-2} & \dots  & d_{1,1} \\
  d_{2,n-1} & d_{2,n-2} & \dots  & d_{2,1} \\
  \vdots &\vdots & \ddots & \vdots \\
  d_{M,n-1} & d_{M,n-2} & \dots  & d_{M,1} \\
\end{bmatrix}  
 \begin{bmatrix}
    f_1(\bm{X}_1) & \dots  & f_M(\bm{X}_1) \\
    f_1(\bm{X}_2) & \dots  & f_M(\bm{X}_2) \\
    \vdots & \ddots & \vdots \\
    f_1(\bm{X}_{n-1}) & \dots  & f_M(\bm{X}_{n-1}) \\       
\end{bmatrix},
\end{displaymath}
and the second term \(\bm{a}_0\bm{F}_n\) is related to the initial approximation of the current step

\begin{displaymath}
\bm{a}_0\bm{F}_n
=  
\begin{bmatrix}
    d_{1,0}\times f_1(\bm{X}^{[\mu-1]}_n)\\
    d_{2,0}\times f_2(\bm{X}^{[\mu-1]}_n)\\
    \vdots \\
    d_{M,0}\times f_M(\bm{X}^{[\mu-1]}_n)
\end{bmatrix}.
\end{displaymath}

Garrappa~\cite{Garrappa2018} suggested the iterative modified Newton-Raphson method for the solution of Eq.~\eqref{eq: IM} with the convergence order \(\mathcal{O}(h^2)\). With having an initial approximation \(\bm{X}_n^{[0]}\) for \(\bm{X}_n\), we can calculate new improved solutions \(\bm{X}_n^{[\mu]}\) by the following iterative formula
\begin{equation}\label{eq: NR}
    \bm{X}_n^{[\mu+1]}=\bm{X}_n^{[\mu]}-\left[ \bm{I}-\bm{a}_0\bm{J_F}(\bm{X}^{[0]}_n)\right] ^{-1} \left( \bm{X}_n^{[\mu]}-\bm{\Psi}_n-\bm{a}_0 \bm{F}_n\right),
\end{equation}
where \(I\) is the identity matrix and \(\bm{J_F}(\bm{X}_n)\) is the Jacobian matrix of \(\bm{F}\) with respect to the variables \(X_i\) at \(t_n\) (or a derivative for one dimensional problems) defined as
\begin{displaymath}
\bm{J_F}(\bm{X}_n)=
\begin{bmatrix}
    \frac{\partial f_1}{\partial X_1} & \frac{\partial f_1}{\partial X_2} \dots  & \frac{\partial f_1}{\partial X_M} \\
    \vdots & \ddots & \vdots \\
    \frac{\partial f_M}{\partial X_1} & \frac{\partial f_M}{\partial X_2} \dots  & \frac{\partial f_M}{\partial X_M} \\
\end{bmatrix}(t_n).
\end{displaymath}
Then, we only need an initial approximation for each step that we achieve from the last evaluated approximation \(\bm{X}_{r+1}^{[0]}={\bf{X}}_r\), (\(r=0,...,n-1\)). It might seem this assumption leads to the method's inefficiency. Still, the high-order convergences of the two blended algorithms and a sufficient number of iterations guarantee efficiency unless the variables \(\bm{X}\) or derivative functions \(\bm{F}\) alter very rapidly (see Example \ref{Sec:LVmulti}) and the step size is not small enough.

\subsection{Fast Fourier Transform}
In both methods, there are convolution sums in the second term of Eq.~\eqref{explicit}, the third term of Eq.~\eqref{implicit}, and the second term of Eq.~\eqref{eq: IM}. Thus, the whole evaluation of the solution requires \(n\) number of operations for grids \(1,...,n\), which proportionally costs \(\mathcal{O}(n^2)\). Hence, the direct calculation of the matrix is not reasonable when the number of grid points \(n\) is sufficiently large. However, Refs.~\cite{Diethelm-FFT} elaborately explained how using FFT can effectively reduce the computational cost proportional to \(\mathcal{O}(n(log_2n)^2)\). 

Here, we briefly describe the exploitation of the FFT algorithm in the PC and NR methods. Suppose 
\begin{math}
Y_n=\sum_{j=1}^n d_{n-j}f_i(\bm{X}_j)
\end{math}
is one of the last sums of our matrix products. The idea is to split each sum into two halves of the computation interval, then Fourier Transform converts the convectional coefficients of one part into frequencies, multiplies once, and finally, converts back the evaluation of the sum. This process requires only \(\mathcal{O}(2n log_2 2n)\) proportional operations, and it can be reduced if we recursively repeated by splitting the interval. However, we need to start the evaluation from the initial time to a smaller length of the interval \(r<n\), 
\begin{displaymath}
Y_q=\sum_{j=1}^q d_{q-j}f_i(\bm{X}_j), \quad q=1,2,...,r. 
\end{displaymath}
This is the initial vector of values that can be directly evaluated, and the length \(r\) is fixed, which could be any small integer number power of two for convenience (where we consider \(r=2^4\)). Then we calculate the values of the following \(r\) points of the interval:
\begin{displaymath}
Y_q=Y_1^r+Y_r^q=\sum_{j=1}^{r} d_{q-j}f_i(\bm{X}_j)+ \sum_{j=r}^{q} d_{q-j}f_i(\bm{X}_j), \quad q=r+1,r+2,...,2r,
\end{displaymath}
where FFT algorithm can evaluate the partial sum \(Y_1^{r}\), with a computational cost proportional to \(\mathcal{O}(2r \log_2 2r)\) instead of  \(\mathcal{O}(r^2)\). For the computation of the next \(2r\) values, we have two splits from $2r$ to $3r$ and from $3r$ to $4r$:
\begin{align*}
    &Y_q= Y_1^{2r}+Y_{2r}^{q} = \sum_{j=1}^{2r} d_{q-j}f_i(\bm{X}_j)+ \sum_{j=2r}^{q} d_{q-j}f_i(\bm{X}_j), \quad q=2r+1,2r+2,...,3r,\\
    & Y_q = Y_1^{2r}+Y_{2r}^{3r}+Y_{3r}^{q} = \sum_{j=1}^{2r} d_{q-j}f_i(\bm{X}_j)+ \sum_{j=2r}^{3r} d_{q-j}f_i(\bm{X}_j) +\sum_{j=3r}^{q} d_{q-j}f_i(\bm{X}_j),\quad q=3r+1,3r+2,...,4r,
\end{align*}
where FFT algorithm can evaluate the partial sums \(Y_1^{2r}\) and \(Y_{2r}^{3r}\), with computational costs proportional to \(\mathcal{O}(4r \log_2 4r)\) and \(\mathcal{O}(2r \log_2 2r)\), respectively. Similarly, this process can be repeated until reaches the \(n\)th point.

\section{Software Implementation with Julia}\label{Sec: software}

We have implemented the numerical methods described in Section \ref{Sec:PC and NR} in the FdeSolver Julia package. It is released with the permissive MIT open-source license, which has been recommended for research software for instance in \cite{morin2012}. The package is also listed on Julia’s General Registry (https://github.com/JuliaRegistries/General). Our current representation is based on the v1.0.7 release of the FdeSolver package. 

The FdeSolver package takes advantage of several features of Julia. This is a compiled language, which means it is interactive and can be used with a read-evaluate-process-loop (REPL) interface. Anyone can contribute to improving the open-source FdeSolver package, and use it in their own applications. Multiple dispatch is one of the useful features that we use in our package. In multiple dispatch a function can have multiple implementations that are allocated for different parameters which would be dispatched at runtime and determined based on the precise parameter types~\cite{Bezanson_Julia}. Multiple dispatch relies on two other performance-enhancing features, composite types, and dynamic types; 1)~There is a unique feature of Julia's composite types (like objects or structs in other languages) in that functions do not get bound to objects nor are they bundled with them. This is essential for multiple dispatch and leads to more flexibility. 2)~Julia can either assign a type to a variable, similar to static programming, or support dynamic types that are determined during execution, contrary to other high-level languages.

\begin{figure}[ht!]
    \centering
    \includegraphics[width=\textwidth]{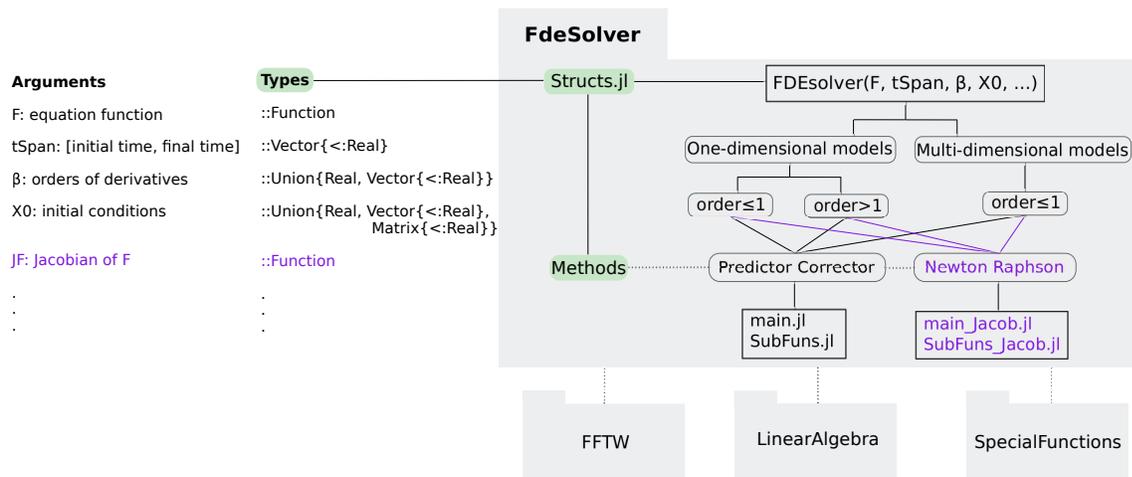}
    \caption{A schematic of FdeSolver package v1.0.7. This package can solve one-dimensional fractional models for the order derivative \(\beta>0\) and multi-dimensional ones for the order derivatives \(0<\beta \leq 1\) of a class of problem \eqref{eq:scalar}. The name of the solver function is FDEsolver. There are ten arguments where the four positional arguments are \textit{F}: the equation function, \textit{tSpan}: the domain of the problem, \(\beta\): the orders of derivatives, and \textit{X0}: the initial conditions and the four optional arguments are \textit{h}: step-size of computation, \textit{nc}: the number of correction for the predictor-corrector method, \textit{itmax}: the maximum number of iterations for the Newton-Raphson method, and \textit{tol}: the tolerance of each iteration or correction. The input \textit{par} is any additional parameters related to the problem and \textit{JF} is a Jacobian function needed for the Netwon-Raphson method. The predictor-corrector method~\ref{Sec: pc} is encoded in main.jl and SubFuns.jl, and the Newton-Raphson method~\ref{Sec: NR} in main\_Jacob.jl and SubFuns\_Jacob.jl.
    Structs.jl leverages multiple dispatch to choose a proper numerical algorithm based on the presence or the lack of a Jacobian function and tune its functionality based on arguments' types. The FFTW, LinearAlgebra, and SpecialFunctions Julia packages support the FdeSolver package.}
    \label{fig:diagram}
\end{figure}

\subsection{Installing FdeSolver}
FdeSolver v1.0.7 is available for Julia version 1 and higher. Julia has a built-in package manager, named ``Pkg'', that can handle operations. By typing a ] in Julia REPL, you can enter the Pkg REPL and type ``add FdeSolver@v1.0.7'', to install the package for this specific version, and use ``up FdeSolver'' or ``rm FdeSolver'' to update the version or remove the package. 

\subsection{Third Party Supporting Packages}
FdeSolver uses FFTW package version 1.2 for fast Fourier transforms, LinearAlgebra for constructing diagonal and identity matrices and applying norm functions, and SpecialFunctions version 1 or 2 for using \(\Gamma\) function.
All dependencies of FdeSolver are listed in Project.toml:
\begin{lstlisting}
[deps]
FFTW = "7a1cc6ca-52ef-59f5-83cd-3a7055c09341"
LinearAlgebra = "37e2e46d-f89d-539d-b4ee-838fcccc9c8e"
SpecialFunctions = "276daf66-3868-5448-9aa4-cd146d93841b"
\end{lstlisting}
and the compatibility constraints for the mentioned dependencies are listed as follows:
\begin{lstlisting}
[compat]
FFTW = "1.2"
SpecialFunctions = "1, 2"
julia = "^1"
\end{lstlisting}

\subsection{Practices for reproducible software development}
FdeSolver follows state-of-the-art methods for building and sharing scientific computing software, as encouraged by the Software Carpentry and Data Carpentry communities~\cite{wilson2017good}. The main repository of the package is created on GitHub (see Sec.~\ref{sec: data-repo}). Unit testing is automatically performed on Ubuntu and MacOS x64 machines via GitHub Actions when the content of one or multiple files is modified. At the same time, Codecov, an open-source third-party service, is used to keep track of the coverage, that is, the percentage of code actually run by the unit tests. Our continuous integration approach simplifies and speeds up the collaborative work around FdeSolver, to which any member of the Julia community can readily contribute. The functionality of our package is thoroughly described in README.md and documented more comprehensively in its official manual, where users can rapidly get started with the four provided usage examples of the solver.

\subsection{FdeSolver Basics}
In Julia, a {\it struct} defines a type that is composed of other types; e.g Float64, Int64, or even other structs. Parametric structs are used throughout FdeSolver to adjust the algorithms for different model classes. 
The types implemented by the package are shown in Fig. \ref{fig:diagram}. The inputs are divided into four sets of arguments:
\begin{itemize}
    \item Positional arguments
    \begin{itemize}
        \item F: derivative function, or the right side of the system of differential equations. Depending on the problem, it could be expressed in the form of a function and return a vector function with the same number of entries of order derivatives. This function can also include a vector of additional parameters.
        \item tSpan: the time along which computation is performed. It must be a vector containing the two values for the initial time and for the final time.
        \item X0: the initial conditions. The values in the type of a row vector for \(\beta \leq 1\) and a matrix for \(\beta> 1\), where each column corresponds to the initial values of one differential equation and each row to its order of derivation.
        \item \(\beta\): the orders derivatives in a type of scalar or vector, where each element corresponds to the order of one differential equation. It could be an integer value.

\begin{lstlisting}
    struct PositionalArguments
       F::Function
       tSpan::Vector{<:Real}
       X0::Union{Real, Vector{<:Real}, Matrix{<:Real}}
       beta::Union{Real, Vector{<:Real}}
    end
\end{lstlisting}
    \end{itemize}
    \item Optional arguments
    \begin{itemize}
        \item h: the step size of the computation. The default is \(2^{-6}\).
        \item nc: the desired number of corrections for the PC method, when there is no Jacobian. The default is 2.
        \item tol: the tolerance of errors taken from the infinity norm of each iteration for the NR method or correction when nc>10 for the PC method. The default is \(10^{-6}\).
        \item itmax: the maximal number of iterations for the NR method, when the user defines a Jacobian. The default is 100.
        
\begin{lstlisting}
    struct OptionalArguments
        h::Float64
        nc::Int64
        tol::Float64
        itmax::Int64
    end
    \end{lstlisting}
    \end{itemize}
    \item JF: the Jacobian of F for switching to the NR method. If it is not provided, the solver will evaluate the solution by the PC method.
    \item par: additional parameters for the functions F and JF.
\end{itemize}

\section{Usage and benchmarks}\label{Sec: usage}

In the following, we present examples using FdeSolver to solve all the model classes presented in Fig. \ref{fig:diagram}.
The user needs to type 
\begin{lstlisting}
using FdeSolver, Plots
\end{lstlisting}
at the top of the code for solving the examples  and plotting the results. Some examples require specific packages that we remark on alongside the related implementations. It is worth noting that we consider the initial conditions of all examples equal to zero, while they could be nonzero as well.

We compare the efficiency (speed and accuracy) of the PC and NR methods of the FdeSolver with four solvers from Matlab counterpart~\cite{Garrappa2018} and six methods provided in another Julia package,
FractionalDiffEq v0.2.11~\cite{SciFracx}, which at the current state provides numerical methods for FDEs but is not applicable for one-dimensional equations. Results suggest that our solver provides a performance that is similar to or better than that of the Matlab counterparts (Figs. \ref{fig:Ex1D}, \ref{fig:ExMD}, and \ref{fig:ExRadom}), and significantly greater accuracy than any available Julia alternatives (Fig. \ref{fig:ExMD}(a)).

Let us define some notations for the methods used for solving the examples. As table \ref{tab:code-List} shows, these notations for the Julia solvers started with the letter J and Matlab solvers include the letter M. J1 indicates the PC method in FdeSlover and it is the same as M1, the method used in Matlab code~\cite{Garrappa2018}. J2 indicates the NR method in FdeSolver, which is the same as M2, the method used in Matlab code. We benchmark these methods for the following five examples besides methods with different algorithms named M3, which is based on the NR method but with PI rectangular rule, and M4, based on the explicit PI rectangular rule without PC~\cite{Garrappa2018}. For Example~\ref{Sec: SIR}, we provide possible comparisons of these methods with six additional methods developed in Julia~\cite{SciFracx}, namely J3, which is based on the same method as J1 and M1, J4, the same method as M4, J5, encoded based on FOTF Toolbox~\cite{FOTF-Toolbox}, J6, J7, and J8, based on fractional linear multistep methods~\cite{GARRAPPA201596}.

\begin{table}[ht!]
{\small{  \caption{A list of solvers used in this study for solving the FDEs. The solvers on the same rows are built on the same algorithms.}
  \label{tab:code-List}
  \begin{tabular}{cccccc}
    \toprule
    notation & FdeSolver & notation & Matlab~\cite{Garrappa2018} & notation & FractionalDiffEq~\cite{SciFracx} \\
    \midrule
    J1 & PC & M1 & FDE\_PI12\_PC.m  &J3 & PIPECE \\
    J2 & NR &  M2 & FDE\_PI2\_Im.m  & - & -\\
    - & - & M3 &FDE\_PI1\_Im.m & - & -\\
    - &- & M4 & FDE\_PI1\_Ex.m & J4 & PIEX \\
    - & -&- & -& J5 & NonLinearAlg\\
    - & -&- & -& J6 & FLMMBDF \\
    - & -&- & -& J7 & FLMMNewtonGregory\\
    - & -&- & -& J8 & FLMMTrap\\
  \bottomrule
\end{tabular}}}
\end{table}

We measure the performance of Matlab's codes by using function {\itshape timeit()} and separately run a group of benchmarks for Julia's codes by using BenchmarkTools.jl package. Benchmark in Julia gives a series of execution times in Nano-seconds. Thus, we take the mean of the series and multiply it by \(10^{-9}\) to fairly compare it with the results from Matlab. A sample of benchmarking M1 and J1 is presented as follows:

\begin{lstlisting}
## julia
    # computing the run time
    t= @benchmark FDEsolver(F, $(tSpan), $(X0), $(beta), $(par), h=$(h)) seconds=1 
    T= mean(t).time / 10^9 # convert from nano seconds to seconds
    #computing the error
    tt, X = FDEsolver(F, tSpan, X0, beta, par, h=h)
    ery=norm(X - map(Exact, tt), 2) # Error: 2-norm
    
\end{lstlisting}
\begin{lstlisting}[language=Matlab]
%% Matlab   
    % computing the run time
    Bench(1,1) = timeit(@() FDE_PI12_PC(beta,F,t0,T,X0,h,param));
    % computing the error
    [t,X]=FDE_PI12_PC(beta,F,t0,T,X0,h,param);
    Bench(1,2)=norm((X-Ex));
\end{lstlisting}
Similarly, we can add other methods and include the settings explained in the examples. 
Finding the exact solution of the multi-dimensional models is difficult or not possible, so we measure the accuracy of the methods by comparing the obtained results with the results secured by fine step size in Matlab: 
\begin{lstlisting}[language=Matlab]
[t,Yex]=FDE_PI2_Im(beta,F,JF,t0,T,x0,2^(-10),par, 1e-12);%Solution with a fine step size in matlab
\end{lstlisting}
and Julia: 
\begin{lstlisting}
t, Yex=FDEsolver(F, tSpan, X0, beta, par, JF = JF, h=2^-10, tol=1e-12)#Solution with a fine step size in julia
\end{lstlisting}
The deviation between the two vector solutions taken from the solvers with a fine step could guarantee or refuse that the approximations are converging to the exact solutions.
 
All the experiments are carried out in Julia Version 1.7.3 and Matlab Version 9.12.0.1975300 (R2022a) Update 3 on a computer equipped with a CPU Intel i7-9750H at 2.60 GHz running under the OS Ubuntu 22.04.1 LTS. Using other versions mentioned above may lead to different results than those presented here. The availability of the source code for replicating the results are presented in Sec.~\ref{sec: data-repo}.

\subsection{One-dimensional models}
We start solving three model examples described by one differential equation. The first two examples have \(\beta \leq 1\) and the third one has \(\beta > 1\).

\subsubsection{Non-stiff example}\label{sec:nonstiff}
Consider the following nonlinear fractional differential equation~\cite{diethelm2002predictor}:

\begin{displaymath}
\mathcal{D}_{t_0}^{\beta}X(t)=\frac{40320}{\Gamma(9-\beta)}t^{8-\beta}-3 \frac{\Gamma(5+\beta/2)}{\Gamma(5-\beta/2)}t^{4-\beta/2}+\frac{9}{4}\Gamma(\beta+1)+\left( \frac{3}{2}t^{\beta/2}-t^{4}\right)^3 -\left(X(t)\right)^{3/2}, \quad 0< \beta \leq 1.
\end{displaymath}

With initial value \(X(0)=0\), the exact solution is 
\begin{math}
X(t)=t^8-3t^{4+\beta/2}+\frac{9}{4}t^\beta.
\end{math}

We need to use SpecialFunctions package for using the gamma function in the equation. Let us suppose the final time is equal to 1 and the order derivative \( \beta=0.5\). The Jacobian of the equation is \begin{math}
J_F=-1.5\sqrt{X}.
\end{math}
Thus, the Julia code for this example is as follows:

\begin{lstlisting}
    using SpecialFunctions
    tSpan = [0, 1]     # [intial time, final time]
    X0 = 0             # intial value
    beta = 0.5            # order of the derivative
    # Equation
    par = beta
    F(t, X, par) = (40320 / gamma(9 - par) * t ^ (8 - par) - 3 * gamma(5 + par / 2)
           / gamma(5 - par / 2) * t ^ (4 - par / 2) + 9/4 * gamma(par + 1) +
           (3 / 2 * t ^ (par / 2) - t ^ 4) ^ 3 - X ^ (3 / 2))
    # Jacobian function
    JF(t, X, par) = -(3 / 2) * X ^ (1 / 2)
    # Solution
    t1, X1 = FDEsolver(F, tSpan, X0, beta, par)
    _, X2= FDEsolver(F, tSpan, X0, beta, par, JF = JF)
\end{lstlisting}
This example has a smooth solution, despite its derivative function displaying a nonlinear and nonsmooth equation. We run the solvers for this range of step size of computations: \(h=2^{-n},\; n=3,...,8\).
Fig. \ref{fig:Ex1D}(a) shows the performance of Matlab and Julia equivalence codes; J1 versus M1 and  J2 versus M2. Hence, in this example, the NR method is more efficient than the PC method, while Julia performs slightly faster than Matlab, and these four solvers execute more accurately than the two other Matlab alternatives, M3 and M4. 

\subsubsection{Stiff example}\label{sec: stiff}
Considering the Jacobian function and using the NR method is recommended for stiff problems~\cite{Garrappa2018} such as:
\begin{displaymath}
\mathcal{D}_{t_0}^{\beta}X(t)=\lambda X(t), \quad X(t_0)=X_0,
\end{displaymath}
where the exact solution is 
\begin{math}
X(t)=E_{\beta}\left((t-t_0)^{\beta}\lambda\right)
\end{math}
in which the Mittag-Leffler function is defined as:
\begin{displaymath}
E_{\beta}(t)=\sum_{k=0}^{\infty}\frac{t^k}{\Gamma(\beta k+1)}.
\end{displaymath}
For this function, we recommend using ``mittleff'' function of  FractionalDiffEq (v0.2.11) package rather than MittagLeffler.jl package since the latter is restricted by compatibility requirements with almost all versions of SpecialFunctions package. 

We solve this problem until \(t=5\) with initial condition \(X(0)=1\), order \(\beta=0.8\), and \(\lambda=-10\). To achieve better results, we set 4 numbers of correction for the PC method.

\begin{lstlisting}
using FractionalDiffEq # to get MittagLeffler function
    [...]
    lambda = -10
    par = lambda
    F(t, X, par)= par * X # equation F
    JF(t, X, par) = par # Jacobian F
    t1, X1 = FDEsolver(F, tSpan, y0, beta, par, nc=4)
    _, X2 = FDEsolver(F, tSpan, y0, beta, par, JF = JF, tol=1e-8)
    # exact solution: mittag-leffler
    Exact(t) = mittleff(beta, lambda * t .^ beta)
    X_exact=map(Exact, t1) # the exact solution corresponding to the time steps t1
\end{lstlisting}
\begin{lstlisting}
using FractionalDiffEq # to get MittagLeffler function
    [...]
    lambda = -10
    par = lambda
    F(t, X, par)= par * X # equation F
    JF(t, X, par) = par # Jacobian F
    t1, X1 = FDEsolver(F, tSpan, y0, beta, par, nc=4)
    _, X2 = FDEsolver(F, tSpan, y0, beta, par, JF = JF, tol=1e-8)
    # exact solution: mittag-leffler
    Exact(t) = mittleff(beta, lambda * t .^ beta)
    X_exact=map(Exact, t1) # the exact solution corresponding to the time steps t1
\end{lstlisting}
This example illustrates the superiority of the NR method. 
We run the solvers for this range of step size of computations: \(h=2^{-n},\; n=3,...,8\).
The PC method diverges for \(h>2^{-4}\), however, the 2-norm error of the NR method for \(h=2^{-3}\) is about 0.16. FdeSolver performs more efficiently and reliably than Matlab counterparts, as Fig. \ref{fig:Ex1D}(b) shows, when the step size is \(h<2^{-7}\) the errors of M1 and M2 increase. 

\subsubsection{High-order example}\label{sec: high-order}
We can solve Example \ref{sec:nonstiff} for an order greater than 1. But, let us consider another popular example, fractional Harmonic motion~\cite{HarmonicMotion}, described as:
\begin{displaymath}
\mathcal{D}_{t_0}^{\beta}X(t)=-\frac{k}{m} X(t), \quad k>0, \quad m>0, \quad 1<\beta \leq 2,
\end{displaymath}
where \(X\) is the displacement from the equilibrium point, \(k\) the spring constant, and \(m\) is the inertial mass of the oscillating body. The exact solution for the \(\beta=2\) is 
\begin{displaymath}
X(t)=X(0)cos(\omega t) + \frac{X'(0)}{\omega}sin(\omega t), \quad \omega=\sqrt{\frac{k}{m}}.
\end{displaymath}
Let us consider \(X(0)=1\) and \(X'(0)=1\) when the parameters are \(k=16\) and \(m=4\), and solve the equation for \(\beta=2\) until \(t=10\). The related code for functions and conditions is represented below.

\begin{figure}[ht!]
    \centering
    \includegraphics[width=.9\textwidth]{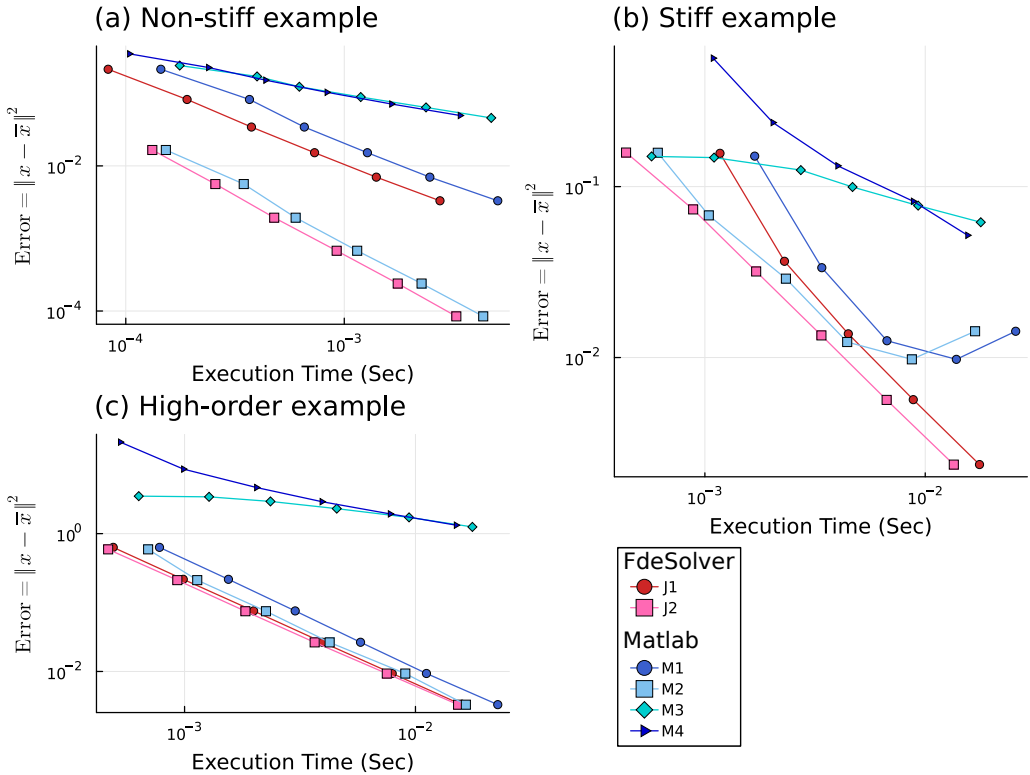}
    \caption{Error versus execution time of FdeSolver and Matlab codes for solving the one-dimensional differential equation of (a) Example~\ref{sec:nonstiff}, (b) Example~\ref{sec: stiff}, and (c) Example~\ref{sec: high-order}. The axes are visualized on a logarithmic scale. The unit of the X-axes is second (Sec) and the Y-axes show the Euclidean norm of the difference between the approximations and exact values. 
    The notations used for the methods of FdeSolver are; J1: the PC method, and J2: the NR method. The notations used for the methods of the Matlab routines are; M1: the PC method, M2: the NR method, M3: the NR method but with the PI rectangular rule, and M4: the explicit PI rectangular rule without PC. 
    }
    \label{fig:Ex1D}
\end{figure}

\begin{lstlisting}
    ## inputs
    tSpan = [0, 10]     # [intial time, final time]
    X0 = [1 1]             # intial value ([of order 0      of order 1])
    beta = 2            # order of the derivative
    par = [16.0, 4.0] # [spring constant for a mass on a spring, inertial mass]
    ## Equations
    function F(t, x, par)
      K, m = par
      - K / m * x
    end
    function JF(t, x, par)
      K, m = par
      - K / m
    end
    [...]
    Yexact = y0[1] .* map(cos, sqrt(par[1] / par[2]) .* t) .+ y0[2] ./ sqrt(par[1] / par[2]) .* map(sin, sqrt(par[1] / par[2]) .* t) # the exact solution map to the corresponding time points
\end{lstlisting}
Fig. \ref{fig:Ex1D}(c) shows the similar performance of FdeSolver and its Matlab counterparts, while both are better than M4 and M3.
We run the solvers for this range of step size of computations: \(h=2^{-n},\; n=2,...,7\).

\subsection{Multi-dimensional systems}\label{sec: mDex}
FdeSolver can solve high dimensional systems of incommensurate fractional differential equations, which means the order derivatives could be unequal. Here, we bring two three-dimensional examples of FDE systems with and without oscillation dynamics to illustrate the performance of the PC and NR methods.

\subsubsection{Non-oscillation example}\label{Sec: SIR}
The susceptible-infected-recovered (SIR) model is the most popular mathematical model for simulating the transmission dynamics of infectious diseases in a population that Kermack and McKendrick have introduced~\cite{Kermack}. Recently, its extensions with fractional orders have been commonly used in different contexts, which motivates us to bring it here as an example and its generalization as an application in Sec. \ref{Sec: Covid}. Suppose the following fractional SIR model~\cite{Saeedian2017PhysRevE} for susceptible \((S)\), infected \((I)\), and recovered \((R)\) individuals with the corresponding order derivatives \(\alpha_S, \alpha_I,\) and \(\alpha_R\):
\begin{flalign*}
    & \mathcal{D}_{t_0}^{\alpha_S} S(t) = -\beta I(t) S(t), \\ 
    & \mathcal{D}_{t_0}^{\alpha_I}I(t) = \beta I(t) S(t) - \gamma I(t), \\ 
    & \mathcal{D}_{t_0}^{\alpha_R}R(t) = \gamma I(t),
\end{flalign*}
where \(\beta\) is the infectious rate and \(\gamma\) the recovery rate. We solve this system until \(t=100\) for initial conditions \(S(0)=1-I(0),\) \(I(0)=0.1\), and \(R(0)=0\), orders \(\alpha_S=0.9, \alpha_I=0.6,\) and \(\alpha_R=0.7\), and parameters \(\beta=0.4\) and \(\gamma=0.04\). Notice that we consider \(tol=10^{-8}\) for the NR method to have a distinguishable result in benchmarking. The related code can be provided below. 
\begin{lstlisting}
## inputs
I0 = 0.1             # intial value of infected
tSpan = [0, 100]       # [intial time, final time]
y0 = [1 - I0, I0, 0]   # initial values [S0,I0,R0]
alpha = [.9, .6, .7]          # order derivatives
## ODE model
par = [0.4, 0.04] # parameters [infectious rate, recovery rate]
function F(t, y, par)
    # parameters
    beta, gamma = par   # infection rate,  recovery rate
    S, I, R = y   # Susceptible, Infectious, Recovered
    # System equation
    dSdt = - beta * S * I
    dIdt = beta * S * I - gamma * I
    dRdt = gamma * I
    return [dSdt, dIdt, dRdt]
end
function JF(t, y, par) # Jacobian of ODE system
    # parameters
    beta, gamma = par   # infection rate,  recovery rate
    S, I, R = y   # Susceptible, Infectious, Recovered
    # System equation
    J11 = - beta * I;    J12 = - beta * S;    J13 =  0
    J21 =  beta * I;    J22 =  beta * S - gamma;    J23 =  0
    J31 =  0;    J32 =  gamma;    J33 =  0
    J = [J11 J12 J13
         J21 J22 J23
         J31 J32 J33]
    return J
end
[...]
\end{lstlisting}

\begin{figure}[ht!]
    \centering
    \includegraphics[width=.9\textwidth]{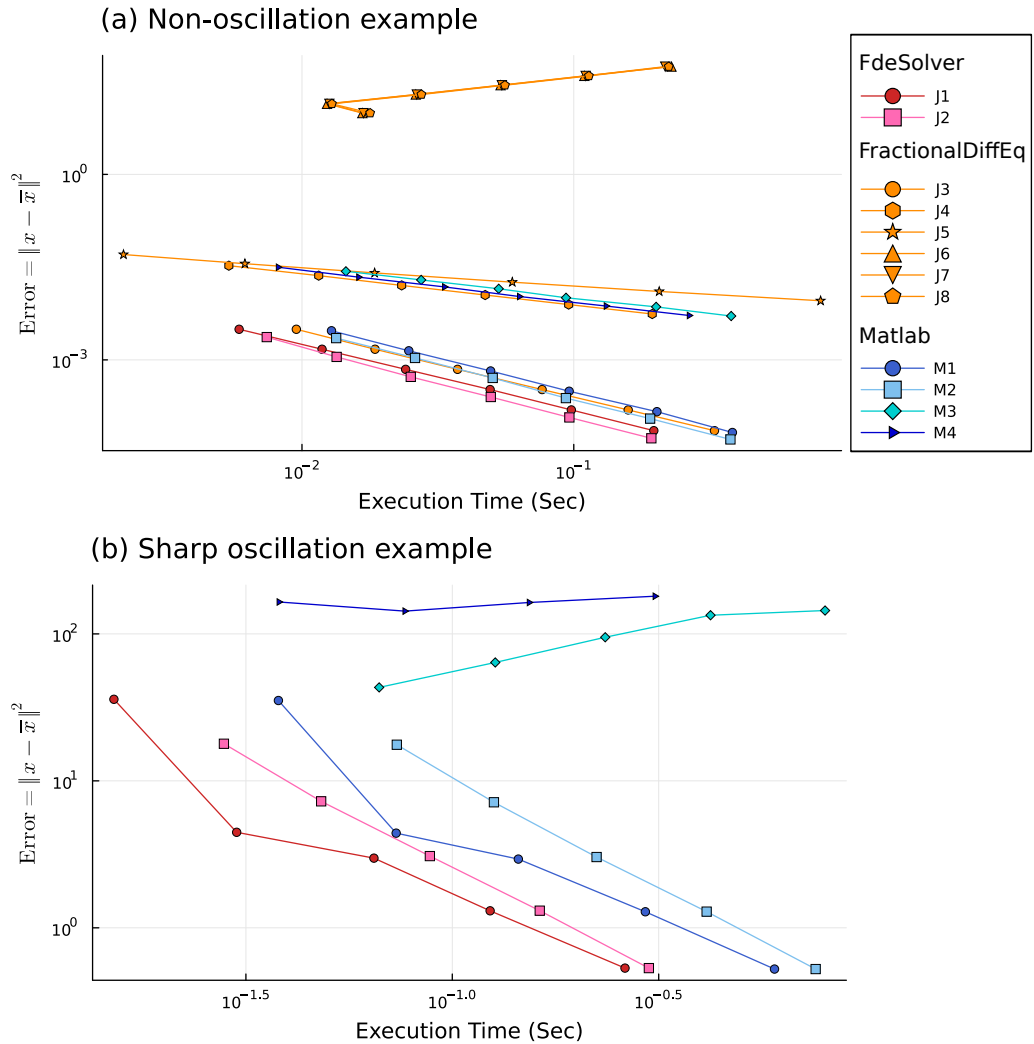}
    \caption{Error versus execution time of Julia and Matlab codes for solving for multi-dimensional systems of (a) Example \ref{Sec: SIR} and (b) Example \ref{Sec:LVmulti}. The axes are visualized on a logarithmic scale. The unit of the X-axes is second (Sec) and the Y-axes show the accuracy by measuring the Euclidean norm of the difference between the approximations and the results with fine step size \(2^{-10}\).
    The notations used for the methods of FdeSolver are; J1: the PC method, and J2: the NR method. The notations used for the methods of FractionalDiffEq are; J3: the PC method, J4: the explicit PI rectangular rule without PC, J5: according to the FOTF Toolbox, J6, J7, and J8: the fractional linear multistep methods. The notations used for the methods of the Matlab routines are; M1: the PC method, M2: the NR method, M3: the NR method but with the PI rectangular rule, and M4: the explicit PI rectangular rule without PC. 
    }
    \label{fig:ExMD}    
\end{figure}

We add six additional methods from FractionalDiffEq to benchmark this example. The methods J6, J7, and J8, from the alternative Julia package, fail, the methods J4 and J5 have a similar performance with M3 and M4, and the method J3 is similar to M1 which is the most reliable solver from FractionalDiffEq (Fig.~\ref{fig:ExMD}(a)).  

Approximations taken from M2 and J2 with fine step size \(h=2^{-10}\) and \(tol=10^{-12}\) are the candidates of the exact solutions for measuring the errors.
The square norm of the difference of these solutions is about \(2.68\times10^{-12}\). We run the solvers for \(h=2^{-n},\; n=2,...,7\) and Fig.~\ref{fig:ExMD}(a) shows FdeSolver superiority. 

\subsubsection{Sharp oscillation example}\label{Sec:LVmulti}
The dynamics of the previous example are smooth, and here we challenge the solver by considering a three-species Lotka-Volterra model~\cite{KhalighiSymmetry} with a sharp oscillation behavior (Fig. \ref{fig:pltDyMD}(b)). This can be defined as:

\begin{align*}
\mathcal{D}_{t_0}^{\beta_1}X_1(t) &= X_1(t)\left( a_1-a_2 X_1(t) - X_2(t) - X_3(t) \right)\\
\mathcal{D}_{t_0}^{\beta_2}X_2(t)&=X_2(t)\left(1- a_3+ a_4X_1(t) \right)\\
\mathcal{D}_{t_0}^{\beta_3}X_3(t)&=X_3(t)\left(1- a_5+a_6 X_1(t) +a_7 X_2(t) \right),
\end{align*}
where \begin{math}
0<\beta_i \leq 1 \; (i=1,2,3),
\end{math}
and all coefficients \begin{math}
a_i>0 \; (i=1,...,7),
\end{math}
and initial values are positive. We solve the system until \(t=60\) for initial conditions \(X_1(0)=X_2(0)=X_3(0)=1\), with parameters \(a_1=a_2=a_3=a_5=a_6=a_7=3\) and \(a_4=5\), and orders \(\beta_1=1,\; \beta_2=0.9\) and \(\beta_3=0.7\). The PC method is more efficient than the NR method for the solution of this example due to fast oscillation dynamics at the initial times, as it is mentioned in Sec. \ref{Sec: NR}.
Hence, we set 4 numbers of corrections for the PC method and \(tol=10^{-8}\) for both methods for the sake of better accuracy and comparison.

This is a challenging example for solvers such that all methods of FractionalDiffEq fail to solve it even for the fine step sizes. Approximations taken from M2 and J2 with fine step size \(h=2^{-10}\) and \(tol=10^{-12}\) are the candidates of the exact solutions for measuring the errors.
The square norm of the difference of these solutions is about \(1.62\times10^{-10}\). We consider \(h=2^{-n},\; n=4,...,8\), and Fig. \ref{fig:ExMD}(b) illustrates that M4 and M3 are not reliable since they are diverging by decreasing the step size. Furthermore, the PC method has a better performance than the NR method in both Julia and Matlab due to the sharp oscillations of the dynamics (see \ref{Sec: NR} and Fig.~\ref{fig:pltDyMD}). 

\subsubsection{Random values}

To assess the generality of the results, let us consider random values for the parameters and conditions of the studied examples. Thus, we run the solvers for the methods J1, J2, M1, M2, M3, and M4 (see Table~\ref{tab:code-List}), to solve Example~\ref{sec:nonstiff} 40 times with random values including $0<\beta \leq 2$ and 30 times each Example of~\ref{sec: stiff}-\ref{Sec: SIR}. We excluded Example~\ref{Sec:LVmulti} due to the difficulty of controlling conditions for making a solvable FDE. Fig.~\ref{fig:ExRadom} shows the distribution of errors and execution times of the solutions of the problems with random conditions, illustrating the similar or improved performance of the FdeSolver Julia package compared to the Matlab alternatives.

\section{Applications}\label{sec: applications}

Let us next present two examples of applications of fractional calculus in community dynamics and epidemiology, where fractional derivatives have been used to describe memory effects, or the influence of past events on population dynamics.~\cite{Saeedian2017PhysRevE,KhalighiPlosCB}. 

\subsection{Simulation of Microbial Community Dynamics}

\begin{figure}[ht!]
    \centering
    \includegraphics[width=\textwidth]{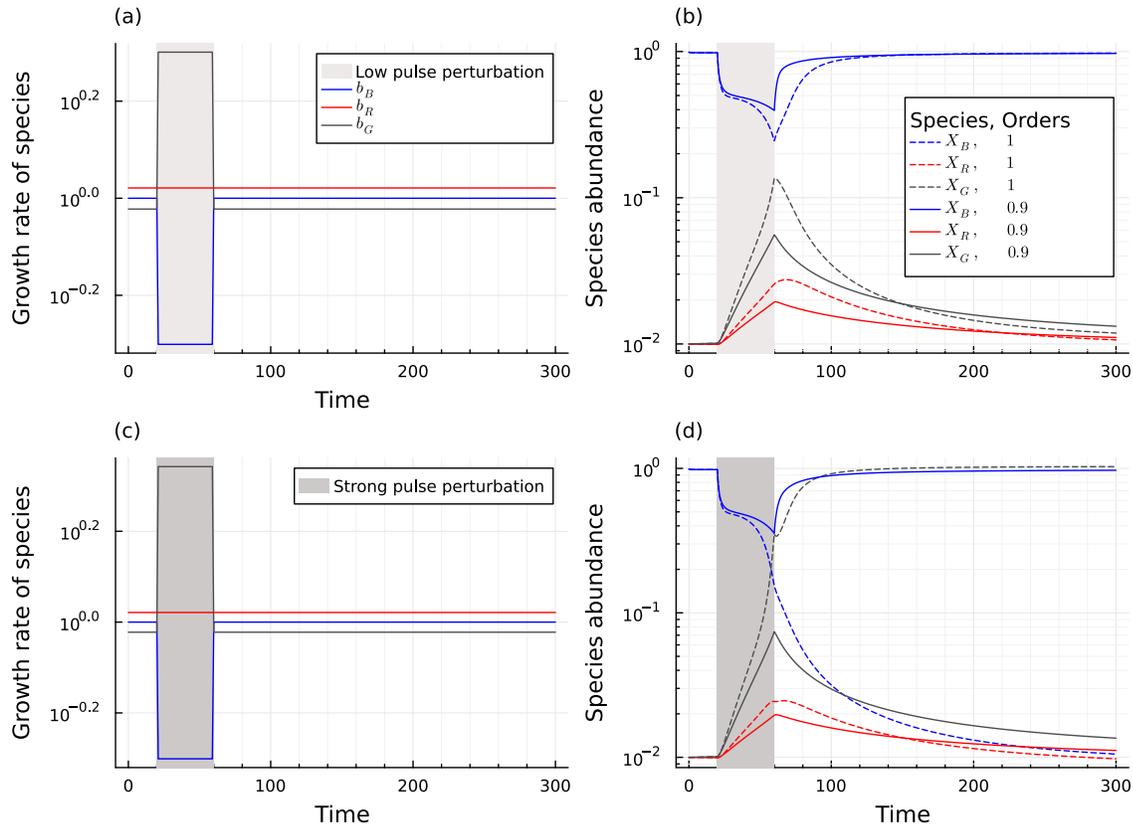}
    \caption{This is the replication of Fig. 2 of Ref.~\cite{KhalighiPlosCB} using the mature FdeSolver Julia package to illustrate the influence of commensurate fractional derivatives on resistance and resilience of the 3-species community model~\eqref{Eq: microb}. Let us suppose \(X_B\), \(X_R\), and \(X_G\) as the abundance of blue, red, and gray species, where indexes \(B, R\), and \(G\) indicate the colors and present through all parameters. We solve the model with initial conditions \(X_B(0)=0.99,\) \(X_R(0)=0.01,\) and \(X_B(0)=0.01\), and parameters \(K_{i,j}=0.1, \forall{i \neq j}\), \(n=2\), \(k_i=1\). The order of the derivatives is set to integer one and fractional values \(\beta_B=\beta_R=\beta_G=0.9\). (a) The growth rates of the system during relaxation are \(b_B=1\), \(b_R=0.95\), and \(b_G=1.05\). A pulse perturbation is applied to the system, during a time window indicated via the grey background, by lowering the growth rate of the blue species (\(b_B=0.5\)) and raising that of the gray species (\(b_G=2\)). (b) The perturbation temporarily moves away the system from the original stable state. However, fractional orders increase resistance to perturbation since the community (solid lines) is not displaced as far from its initial state compared to the dynamics of the system with integer order (dashed lines). (c) A slightly stronger pulse perturbation is applied to the growth rate of the gray species (\(b_G=2.2\)). (d) It triggers a shift toward an alternative stable state dominated by the gray species in the system with integer order (dashed lines). However, fractional orders can also entirely prevent a state shift (solid lines) by increasing both resistance and resilience to perturbation.}
    \label{fig:microb}
\end{figure}

In microbial ecology, long-term memory has been observed in the context of antibiotic-tolerant persister cells~\cite{simsek2019PNAS}, which results from phenotypic heterogeneity~\cite{Gokhale2021PlosCB, Power2015} that has a power-law scale feature. In our recent study, we numerically investigated the effect of memory in ecosystem dynamics of interacting communities~\cite{KhalighiPlosCB}, using some of the methods that ultimately converged in the FdeSolver package. We investigated extensions to interaction models including the generalized Lotka-Volterra model with fractional orders described as:

\begin{equation}\label{Eq: microb}
    \mathcal{D}_{t_0}^{\beta_i}=X_i\left(b_i f_i (\{X_k\})-k_iX_i\right), \quad f_i(\{X_k\})=\prod\limits_{\substack{k=1 \\ k \neq i}}^{N} \frac{K_{ik}^{n}}{K_{ik}^{n}+X_{k}^{n}}, \quad i=1,...,N,
\end{equation}

where \(N\) is the number of species, \(X_i\) species abundances, \(b_i\) growth rates,  \(k_i\) death rates, and \(f_i\) inhibition functions such that \(K_{ij}\) and \(n\) denote interaction constants and Hill coefficients~\cite{gonze2017multi}. If we determine \(1-\beta_i\) as the memory index of the dynamics of species abundance \(i\), then a decrease in the order of the corresponding derivative leads to an increase in its memory.

The FdeSolver package thus allows us to replicate the simulation results that we recently reported~\cite{KhalighiPlosCB}, to  demonstrate and assess system dynamics under pulse and periodic perturbations in the presence of ecological memory (see Fig.~\ref{fig:microb}).

\subsection{Epidemiological analysis of Covid-19 transmission dynamics}\label{Sec: Covid}

Following the global spread of the Coronavirus pandemic in 2019, infectious disease transmission analysis has become very popular, and a large amount of time series data on Covid-19 occurrences has been published. Although this stimulated progress in epidemiological models, predicting the spread of the pandemic proved to be challenging. We demonstrate how the FdeSolver package can speed up simulations and improve model fits in real data. 

Methods for modelling general transmission dynamics are available; for instance, the Julia package Pathogen.jl~\cite{Pathogen} provides many tools to simulate and infer transmission networks and to model parameters of infectious disease spread. Our work can complement and extend such models by providing new tools for incorporating memory effects. These effects have been linked to the spread and control of epidemics, for instance, by changes in precautionary measures, such as vaccinations and lockdowns~\cite{Saeedian2017PhysRevE}. The Covid-19 evolution was suggested to exhibit power-law scaling features~\cite{JAHANSHAHI}, which motivates the use of Caputo derivatives.

A recent compartmental model with super-spreader class~\cite{Covid_Model_Main} was extended with fractional order derivatives~\cite{Covid-Model-Fractional, Covid-Model-Frac-Analysis}. In this model, the population is constant and is divided into eight epidemiological compartments: susceptible individuals (\(S\)), exposed individuals (\(E\)), symptomatic and infectious individuals (\(I\)), super-spreaders individuals (\(P\)), infectious but asymptomatic individuals (\(A\)), hospitalized individuals (\(H\)), recovery individuals (\(R\)), and dead individuals (\(F\)). Here, we investigate an incommensurate fractional orders form of the model defined as:

\begin{align}\label{eq: Covid}
\begin{split}
    & \mathcal{D}_{t_0}^{\alpha_S}S(t) = -\beta \frac{I}{N}S-l \beta \frac{H}{N} S -\beta' \frac{P}{N} S,\\
    & \mathcal{D}_{t_0}^{\alpha_E}E(t) = \beta \frac{I}{N}S + l \beta \frac{H}{N} S +\beta' \frac{P}{N} S -\kappa E,\\
    & \mathcal{D}_{t_0}^{\alpha_I}I(t) = \kappa \rho_1 E - (\gamma_a + \gamma_i)I -\delta_i I,\\
    & \mathcal{D}_{t_0}^{\alpha_P}P(t) = \kappa \rho_2 E - (\gamma_a + \gamma_i)P -\delta_p P,\\
    & \mathcal{D}_{t_0}^{\alpha_A}A(t) = \kappa (1- \rho_1 -\rho_2)E,\\
    & \mathcal{D}_{t_0}^{\alpha_H}H(t) = \gamma_a (I+P) - \gamma_r H - \delta_h H,\\
    & \mathcal{D}_{t_0}^{\alpha_R}R(t) = \gamma_i (I+P) + \gamma_r H,\\
    & \mathcal{D}_{t_0}^{\alpha_F}F(t) = \delta_i I + \delta_p P + \delta_h H,\\
\end{split}
\end{align}
where \(0<\alpha_i \leq 1\) (\(i=S, E, I, P, A, H, R, F\)) are the order derivatives corresponding to the compartments and not necessarily equal, and \(N\) is total population, where \(N = S + E + I + P + A + H + R + F\). The parameters with their values are described in Table~\ref{tab:parCovid}.
\begin{table}[ht!]
{\small{  \caption{Description of the parameters used in model~\eqref{eq: Covid} and their values taken from Refs.~\cite{Covid-Model-Frac-Analysis,Covid-Model-Fractional}.}
  \label{tab:parCovid}
  \begin{tabular}{llll}
    \toprule
    Notation & Description & Value & Units\\
    \midrule
     \(\beta\) & Transmission coeﬃcient from infected individuals & fitted & day\(^{-1}\)\\
    \(\kappa\) & Rate at which exposed become infectious & 0.25 & day\(^{-1}\)\\
    \(l\) & Relative transmissibility of hospitalized patients & 1.56 & dimensionless\\
    \(\beta'\) & High transmission coeﬃcient due to superspreaders & 7.65 & day\(^{-1}\)\\
    \(\kappa\) & Rate at which exposed become infectious &2.55 & day\(^{-1}\)\\
    \(\rho_1\) & Rate at which exposed people become infected \(I\) & 0.58 & dimensionless\\ 
    \(\rho_2\) & Rate at which exposed people become super-spreaders & 0.001 & dimensionless\\ 
    \(\gamma_a\) & Rate of being hospitalized & 0.94 & day\(^{-1}\)\\
    \(\gamma_i\) & Recovery rate without being hospitalized & 0.27 & day\(^{-1}\)\\ 
    \(\gamma_r\) & Recovery rate of hospitalized patients & 0.50 & day\(^{-1}\)\\ 
    \(\delta_i\) & Disease induced death rate due to infected class & 1/23 & day\(^{-1}\)\\ 
    \(\delta_p\) & Disease induced death rate due to super-spreaders & 1/23 & day\(^{-1}\)\\
    \(\delta_h\) & Disease induced death rate due to hospitalized class & 1/23 & day\(^{-1}\)\\
  \bottomrule
\end{tabular}}}
\end{table}

\begin{figure}[ht!]
    \centering
    \includegraphics[width=.9\textwidth]{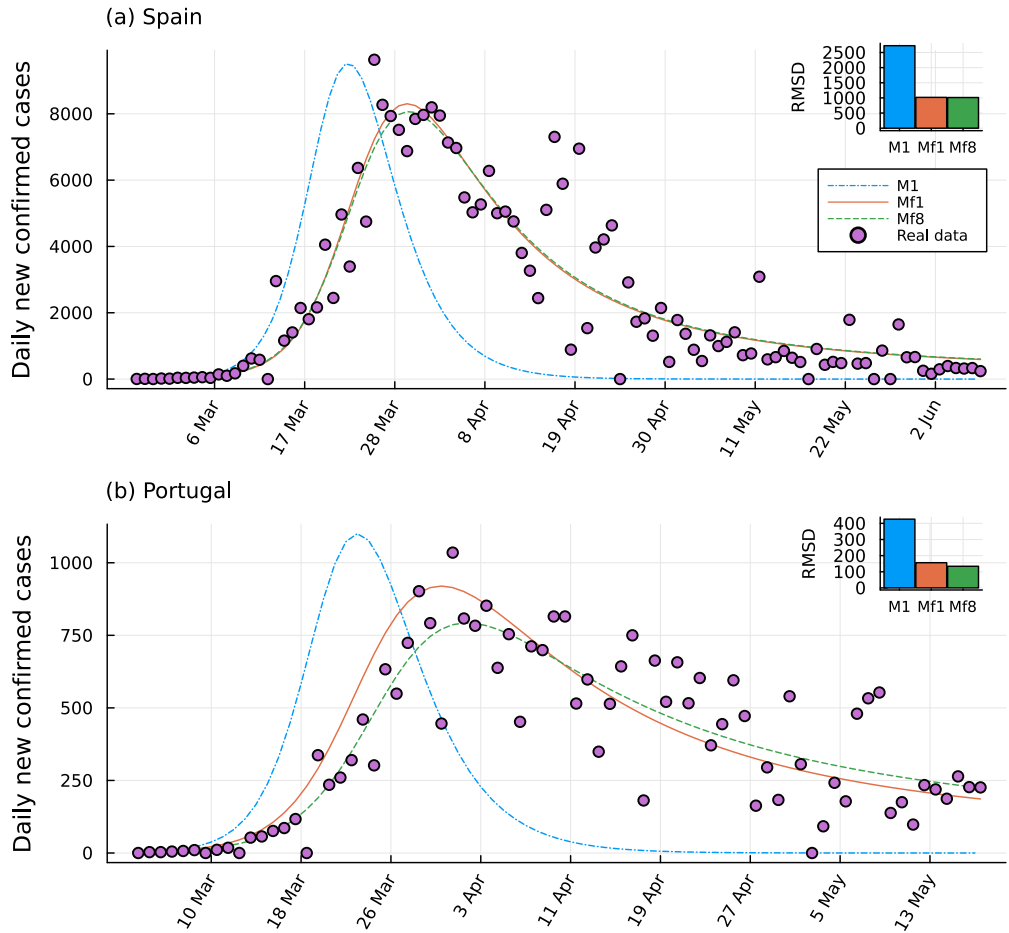}
    \caption{The comparison between real data on the daily new confirmed cases as retrieved from CSSE~\cite{Data-Covid} and the estimation of \(I+P+H\) from the system equation \eqref{eq: Covid} with integer order derivatives (model M1 shown by the dash-dotted blue line),  commensurate fractional  derivatives (model Mf1 shown by the solid orange line), and incommensurate fractional  derivatives (model Mf8 shown by the green dashed line). The circles indicate the real data for (a) Spain from 25 February to 7 June and (b) Portugal from 3 March to 17 May. The infection rate \(\beta\) is fitted for the M1 model, the commensurate order for the Mf1 model together with \(\beta\) are optimized, and finally, the eight incommensurate orders and \(\beta\) for model Mf8 are optimized to fit real data. The specified error bars are based on root mean square deviation (RMSD) such that the fewer values indicate the less residual variance. The fitted values are listed in Table~\ref{tab:fitValues}.}
    \label{fig:Spn_Prt}
\end{figure}

Let us consider the time series of the daily new conﬁrmed cases of Spain and Portugal from cumulative conﬁrmed cases reported by the Center for Systems Science and Engineering (CSSE) at Johns Hopkins University~\cite{Data-Covid}. We take into account three types of the system equation \eqref{eq: Covid} in terms of the derivatives: with integer orders (M1), commensurate fractional orders (Mf1), and incommensurate fractional orders (Mf8). The approximated values for the parameters, the initial conditions, and the population size are taken from Refs.~\cite{Covid-Model-Frac-Analysis,Covid-Model-Fractional}. 

We estimate the parameter infection rate \(\beta\) by fitting the models to the data. For the fractional models (Mf1 and Mf8) we optimize the values of order derivatives together with the parameter \(\beta\) to achieve the best fit. All the fitted values are listed in Table~\ref{tab:fitValues}. The residual of the fitting is measured by root mean square deviation (RMSD) defined as:
\begin{math}
\rm{RMSD}(Y,\widehat{Y})=\sqrt{\frac{1}{n}\sum_{t=1}^n(\widehat{Y_t}-Y_t)^2},
\end{math} 
where \(n\) is the number of data points, \(Y\) the vector of approximations, and \(\widehat{Y}\) real data. 

Fig.~\ref{fig:Spn_Prt} illustrates the simulations of the models along with the real data. The residual errors imply that fitting one parameter of model M1 is not enough, whereas additionally optimizing the order derivatives can remarkably improve the flexibility and, consequently, the accuracy of the model. 

We use the StatsBase.jl package to apply the RMSD for the real values and the sum of all infectious compartments \(I + P + H\) approximated from the models. To minimize the residual, we use two functions from the Optim.jl package~\cite{optim}: (L)BFGS, which is based on the Broyden–Fletcher–Goldfarb–Shanno algorithm, and SAMIN, which is based on Simulated Annealing for problems with bounds constraints.  

\begin{table}[ht!]
{\small{  \caption{Optimized values of the order derivatives \(\alpha_{S,E,I,P,A,H,R,F}\) and parameter \(\beta\) of system equation \eqref{eq: Covid} with integer orders (model M1) commensurate orders (model Mf1) incommensurate orders (model Mf8) for fitting Covid-19 data from Spain and Portugal with their errors measured by root mean square deviation (RMSD).}
  \label{tab:fitValues}
  \begin{tabular}{ccccccc}
    \toprule
     & \multicolumn{3}{c}{Fitted values for Spain's data}&\multicolumn{3}{c}{Fitted values for Portugal's data}\\
    Models & M1 & Mf1 & Mf8 & M1 & Mf1 & Mf8\\
    \midrule
    \(\beta\)& 1.852210 & 2.503197 & 2.520522 & 1.638658 & 2.548180 & 2.841678\\
    \(\alpha_{S,E,I,P,A,H,R,F}\)& -&  0.829366&- & -& 0.775394 &-\\
    \(\alpha_{S}\)& -& - & 0.829128 &- & - & 0.764514 \\
    \(\alpha_{E}\)&- & - & 0.865748 &- & - & 0.724861 \\
    \(\alpha_{I}\)& -& - & 0.685968 &- & - & 0.875856 \\
    \(\alpha_{P}\)& -&-  & 0.500000 & -& - & 1.000000 \\
    \(\alpha_{A}\)& -&  -& 0.749881 & -&  -& 0.749835 \\
    \(\alpha_{H}\)& -& - & 0.809659 & -& - & 0.611715 \\
    \(\alpha_{R}\)& -& - & 0.749881 & -&  -& 0.749835 \\
    \(\alpha_{F}\)&- &  -& 0.749881 & -& - & 0.749835 \\
   \midrule
  RMSD &0.282941 &  0.105667& 0.105163 & 0.411630 & 0.150977 & 0.129533 \\
  \bottomrule
\end{tabular}}}
\end{table}

\section{Conclusion}\label{sec: Conclusion}

Efficient methods for fitting fractional differential equation models can provide valuable tools for dynamical systems analysis in many application fields. The limited availability of accessible tools for fractional calculus in modern open-source computational languages, such as Julia, has slowed down the adoption and development of these methods as part of the broader data science ecosystem. We have introduced a new Julia package that helps to fill this gap through efficient, adaptable, and integrated utilities that facilitate numerical analysis and simulation of FDE models.

Our newly introduced FdeSolver package provides numerical solutions for FDEs with Caputo derivatives in univariate and multivariate systems. The algorithm is based on PI rules, and we have improved its computational efficiency by incorporating the FFT technique. The benchmarking experiments showed a better performance in speed and accuracy as compared with the currently available alternatives in Matlab and Julia. Moreover, we have demonstrated two applications of complex systems for simulating and fitting the models to real data. The package adheres to good practices in open research software development, including open licensing, unit testing, and a comprehensive documentation with reproducible examples.

Future extensions could include fractional linear multistep methods, whose efficiency has been shown to outperform PI rules~\cite{GARRAPPA201596}. We are also seeking to add support for models with time-varying derivatives, thus enhancing the overall applicability across a broader range of problems, e.g. Caputo fractional SIR models with multiple memory effects for fitting Covid-19 transmission~\cite{JAHANSHAHI}. The fractional models are more flexible than their integer-order counterparts, which can potentially lead to either perfect fitting or overfitting. Therefore, automated cross-validation and other means to quantify and control overfitting in parameter optimization will need to be further developed. Moreover, taking advantage of the Bayesian approximation from the Turing.jl~\cite{Turing} package is another promising feature for modelling the uncertainty of order derivatives. Regarding package performance, adding GPU support~\cite{Cuda} could further help to accelerate the speed of computation. Such future developments will further integrate FdeSolver with other Julia packages and deliver its advantageous functionality to the wide Julia ecosystem. 

\section{Code and data Availability}\label{sec: data-repo}
The FdeSolver package is available on GitHub, and accessible via the permanent Zenodo DOI: \url{https://doi.org/10.5281/zenodo.7462094}, and all computational results including all data and code used for running the simulations and generating the figures are accessible from Zenodo DOI: \url{https://doi.org/10.5281/zenodo.7473300}.

\begin{acks}
The research was supported by the Academy of Finland (URL: https://www.aka.fi; decision 330887 to LL, MK), the University of Turku Graduate School Doctoral Programme in Technology (UTUGS/DPT)\\(URL: https://www.utu.fi/en/research/utugs/dpt; to MK), and the Baltic Science Network Mobility Programme for Research Internships (BARI) (URL: https://www.baltic-science.org; to GB).
\end{acks}



\bibliographystyle{ACM-Reference-Format}
\bibliography{ref}

\appendix

\section{Supplementary figures}
There are two additional figures in this Appendix, one showing the dynamics of the multi-dimensional example systems (Fig.~\ref{fig:pltDyMD}) and the other showing the speed and accuracy of solving the problems with random conditions (Fig.~\ref{fig:ExRadom}).

\renewcommand{\thefigure}{A\arabic{figure}}
\setcounter{figure}{0}

\begin{figure}[th!]
    \centering
    \includegraphics[width=\textwidth]{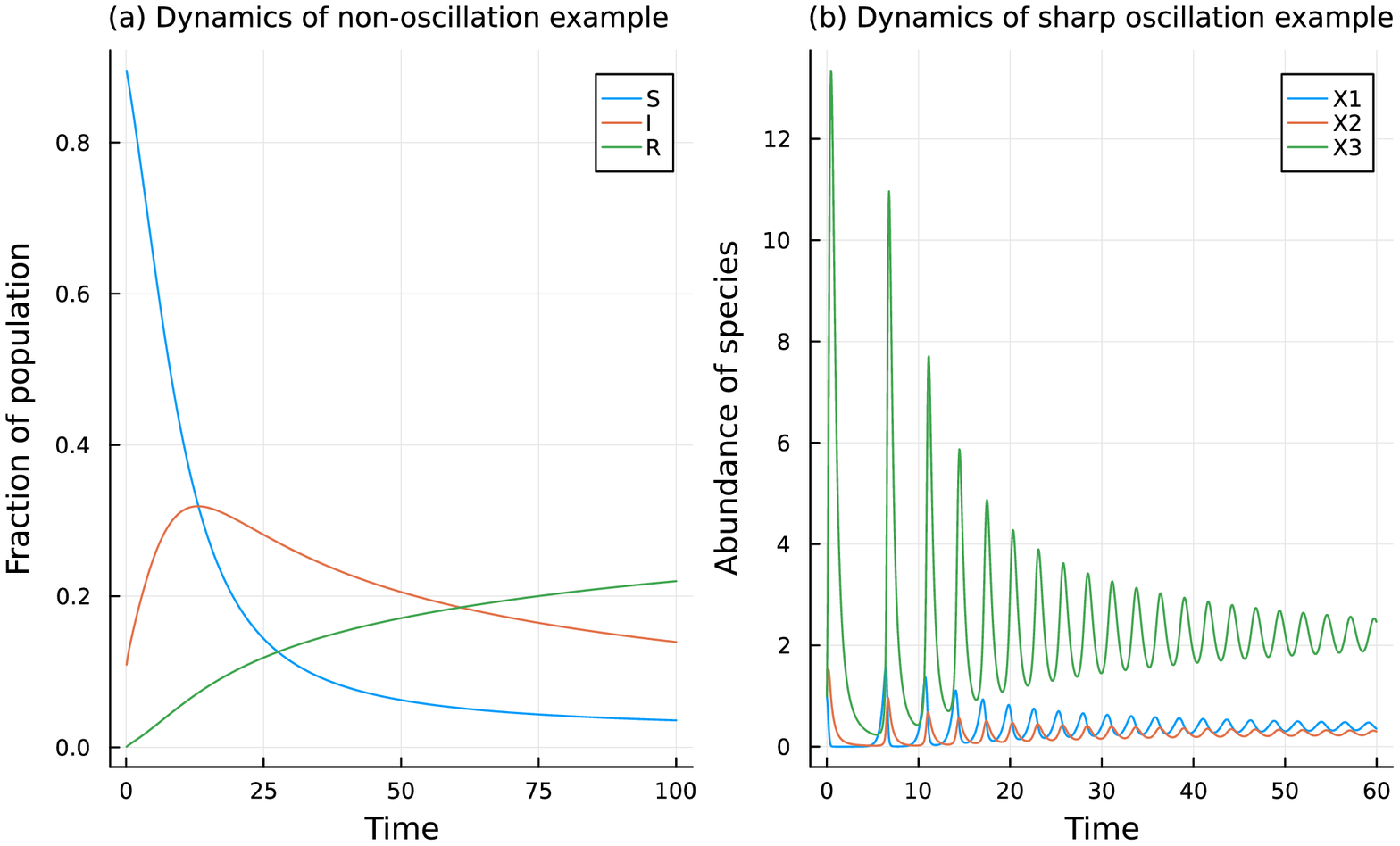}
    \caption{Illustration of the dynamics of multi-dimensional examples. (a) Non-oscillation dynamics of the susceptible \((S)\), infected \((I)\), and recovered \((R)\) population fractions for Example \ref{Sec: SIR}. (b) Sharp oscillation dynamics of the abundance of species \(X_1\), \(X_2\), and \(X_3\) for Example \ref{Sec:LVmulti}.}
    \label{fig:pltDyMD}
\end{figure}

\begin{figure}[th!]
    \centering
    \includegraphics[width=\textwidth]{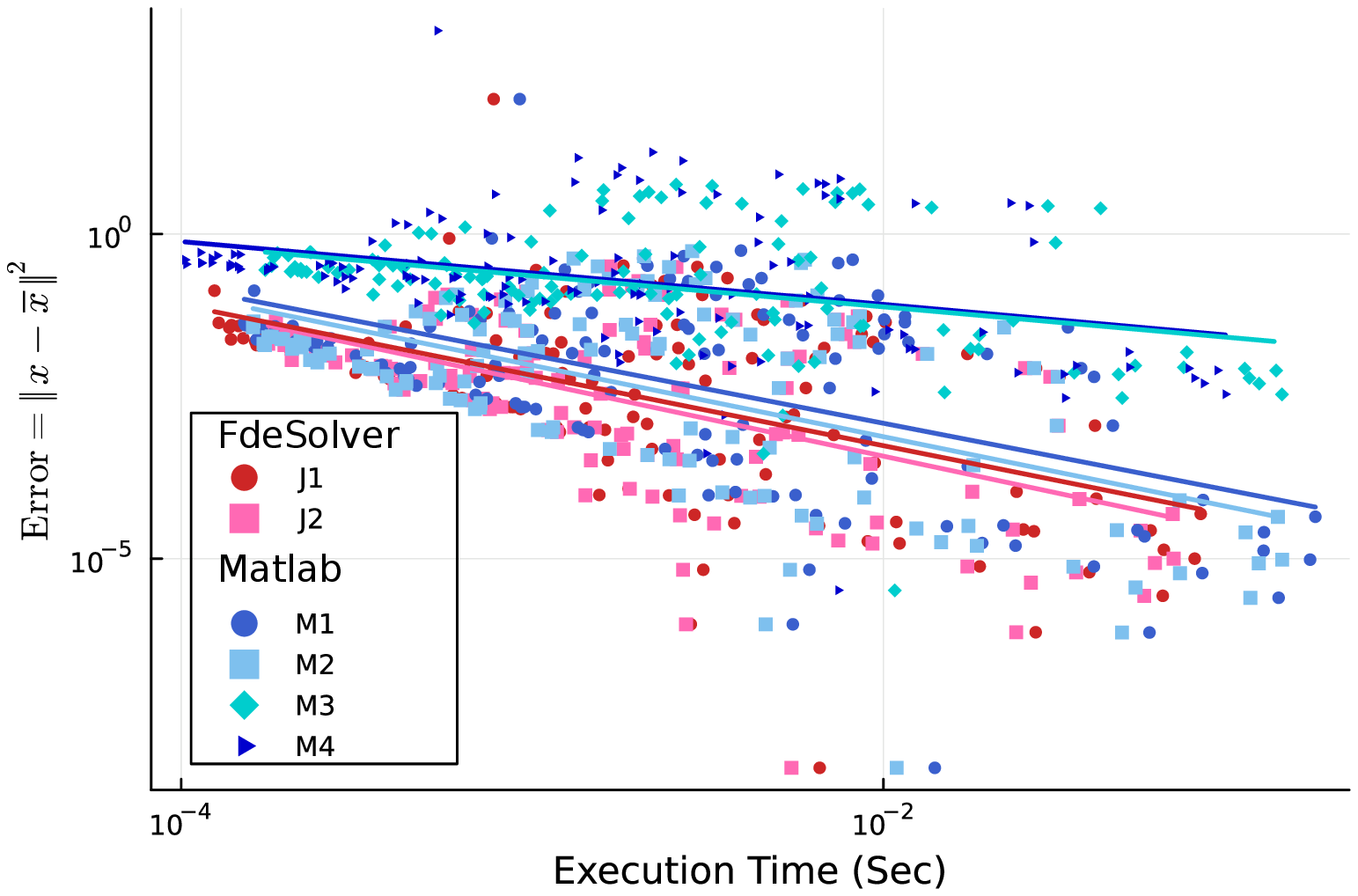}
    \caption{Scatter plot of error versus execution time for several solutions of Examples \ref{sec:nonstiff}-\ref{Sec: SIR} with random parameters and conditions. The axes are visualized on a logarithmic scale. The unit of the X-axes is second (Sec) and the Y-axes show the Euclidean norm of the difference between the approximations and exact values. The lines on the background are the linear fit of the logarithm of the scatter points and the color of the lines corresponds to the methods with the same color.
    The notations used for the methods of FdeSolver are; J1: the PC method, and J2: the NR method. The notations used for the methods of the Matlab routines are; M1: the PC method, M2: the NR method, M3: the NR method but with the PI rectangular rule, and M4: the explicit PI rectangular rule without PC. 
    }
    \label{fig:ExRadom}
\end{figure}

\end{document}